\documentclass[a4paper,12pt]{article}

\usepackage{float}
\usepackage{amsmath} 
\usepackage{amssymb}
\usepackage{amsfonts}
\usepackage{epstopdf}
\usepackage{color,fancybox,graphicx}
\usepackage{url}
\usepackage{psfrag}
\usepackage{color}
\usepackage{pifont}

\usepackage{natbib}
\bibpunct{(}{)}{;}{a}{,}{,}
\usepackage{hyperref}
\hypersetup{
  colorlinks = true,
  urlcolor = blue, 
  linkcolor = blue,
  citecolor = blue,
}

\numberwithin{equation}{section}

\newtheorem{rem}{Remark}[section]

\newtheorem{lem}{Lemma}[section]

\newtheorem{pro}{Proposition}[section]
\newtheorem{theo}{Theorem}[section]

\newcommand{\N}{\mathbb N}
\newcommand{\R}{\mathbb R}
\newcommand{\prob}{\mathbb{P}}
\newcommand{\esp}{\mathbb{E}}

\newcommand{\ra}{\rightarrow}
\newcommand{\e}{\varepsilon}

\definecolor{gris25}{gray}{0.90}

\begin{document}

\begin{center}

{\Large 
\textbf{\textsf{Cox Process Functional Learning}}}
\medskip

\medskip

\end{center}
{\bf G\'erard Biau\\
{\it Universit\'e Pierre et Marie Curie\footnote{Research partially supported by the French National Research Agency (grant ANR-09-BLAN-0051-02 ``CLARA'') and by the Institut universitaire de France.} \& Ecole Normale Sup{\'e}rieure\footnote{Research carried out within the INRIA project ``CLASSIC'' hosted by Ecole Normale Sup{\'e}rieure and CNRS.}, France}}\\
gerard.biau@upmc.fr
\bigskip

{\bf Beno\^{\i}t Cadre}\\
{\it IRMAR, ENS Rennes, CNRS, UEB, France\footnote{Research sponsored by the French National Research Agency (grant ANR-09-BLAN-0051-02 ``CLARA'').}}\\
benoit.cadre@ens-rennes.fr
\bigskip

{\bf Quentin Paris}\\
{\it CREST, ENSAE-ParisTech, France}\\
quentin.paris@ensae.fr
\medskip

\begin{abstract}
\noindent  {\rm This article addresses the problem of functional supervised classification of Cox process trajectories, whose random intensity is driven by some exogenous random covariable. The classification task is achieved through a regularized convex empirical risk minimization procedure, and a nonasymptotic oracle inequality is derived. We show that the algorithm provides a Bayes-risk consistent classifier. Furthermore, it is proved that the classifier converges at a rate which adapts to the unknown regularity of the intensity process. Our results are obtained by taking advantage of martingale and stochastic calculus arguments, which are natural in this context and fully exploit the functional nature of the problem. 
\medskip

\medskip

}

\end{abstract}

\section{Introduction}
\subsection{Functional classification and Cox processes}
In supervised classification one considers a random pair $(X,Y)$, where $X$ takes values in some space $\mathcal X$ and $Y$ takes only finitely values, say -1 or 1 to simplify. Given a learning sample $\mathcal D_n=\{(X_1,Y_1), \hdots, (X_n,Y_n)\}$ of i.i.d.~copies of $(X,Y)$ observed in the past, the aim is to predict the value of $Y$ associated with a new value of $X$. In medicine, for example, we specifically want to evaluate patients according to their disease risk, and the typical questions for classification are: ``Is this person affected?'', ``Will this patient respond to the treatment?'', or ``Will this patient have serious side effects from using the drug?''---in all these cases, a yes/no or $-1/1$ decision has to be made.
\medskip

The classification task is generally achieved by designing a decision rule (also called classifier) $g_n:\mathcal X\to \{-1,1\}$, which represents our guess on the label $Y$ of $X$ (the subscript $n$ in $g_n$ means that the classifier measurably depends upon the sample). As the pair $(X,Y)$ is random, an error occurs whenever $g_n(X)$ differs from $Y$, and the probability of error of the rule $g_n$ is
$$L(g_n)=\mathbb P\left(g_n(X)\ne Y|\mathcal D_n\right).$$
The Bayes rule $g^\star$, defined by 
\begin{equation*}\label{bayesrule}
g^\star(x)= \left\{\begin{array}{ll}
1 \mbox{ if } \mathbb P(Y=1|X=x)\ge \mathbb P(Y=-1|X=x)\\
0 \mbox{ otherwise},
\end{array}\right.
\end{equation*}
has the smallest probability of error, in the sense that $L(g^\star)\le L(g)$ for any classifier $g$ \citep[see, e.g.,][]{DeGyLu96}.
\medskip

In the classical statistical setting, each observation $X_i$ is a collection of numerical measurements represented by a
$d$-dimensional vector. However, in an increasing number of application domains, input
data are in the form of random functions rather than standard vectors, thereby turning the classification task into a functional data analysis problem. Here, the vocable ``random functions'' means that the variables $X_i$'s take values in a space $\mathcal X$ of functions rather than $\mathbb R^d$, equipped with an appropriate topology. Thus, in this context, the challenge is to design classification rules which exploit the functional nature of the $X_i$'s, and this calls for new methodological concepts. Accordingly, the last few years have witnessed important developments in both the theory and practice of functional
data analysis, and numerous procedures have been adapted to handle functional inputs. The books by \citet{RaSi02,RaSi05} and \citet{Fe11} provide a presentation of the area, and the survey of \citet{BaCuFr11} offers some essential references for functional supervised classification.
\medskip

Curiously, despite a huge research activity in the field, few attempts have been made to connect the area of functional data analysis  with the theory of stochastic processes, which also deals with the analysis of time-dependent quantities \cite[interesting ideas towards this direction are included in][]{IlBeCrSt06,BaCuCu11,Ca13,ShMuZh13}. As advocated in the present paper, stochastic calculus theory can be used efficiently to analyse Cox models and may serve as a starting point for more exchanges between the two fields.
\medskip 

To motivate the use of Cox models for classification, consider for instance a sample of AIDS patients observed until time $T$. Assume that, for each of them, we know the dates of visits to the hospital, a bunch of personal data (such as gender, distance from home to hospital, etc.), as well as the diagnostic (-1=aggravation, 1=remission, for example). Based on this learning sample, a classification strategy aims at predicting the $\pm 1$ evolution diagnostic of a new patient. In this time-dependent setting, $\mathcal X$ is the set of counting paths on $[0,T]$ (that is, right-continuous and piecewise constant paths on $[0,T]$ starting at 0, and with jump size 1), and a relevant model for $X$ is a mixture of two Cox processes (or doubly stochastic Poisson processes) with (random) intensities $\lambda_+=(\lambda_{+,t})_{t\in [0,T]}$ and $\lambda_-=(\lambda_{-,t})_{t\in [0,T]}$. In other words, conditionally on $Y=1$ (resp., $Y=-1$), the law of $X$ given $\lambda_+$ (resp., $\lambda_-$) is the law of a Poisson process with intensity $\lambda_+$ (resp., $\lambda_-$). (\citealp[For more information on Cox processes, we refer the reader to the original paper by][]{Co55}; see also the book by \citealp[][for an overview of the application areas of these processes.]{BeKo02}) Compared to a Poisson process, the benefit of the random intensity lies in the fact that the statistician can take into account the auxiliary information carried by the personal data of the patients. \medskip

As we shall see, because of a martingale property of Cox processes, stochastic calculus proves to be a natural and efficient tool to investigate this classification problem.  It is stressed that the originality of our work is that it takes advantage of the theory of stochastic processes to handle a functional data analysis problem---in that sense, it differs from other studies devoted to nonparametric estimation of Cox process intensity \citep[see for instance][and the references therein]{HaRBRi13}.

\subsection{Classification strategy} 
In the sequel, $T>0$ is fixed and $\mathcal X$ stands for the set of counting paths on $[0,T]$. We consider a prototype random triplet $(X,Z,Y)$, where $Y$ is a binary label taking 
the values $\pm1$ with respective positive probabilities $p_+$ and $p_-$  ($p_++p_-=1$). In this model, $Z=(Z_t)_{t\in [0,T]}$ plays the role of a $d$-dimensional random covariable (process), whereas $X=(X_t)_{t\in [0,T]}$ is a mixture of two Cox processes, both being adapted with respect to the same filtration. More specifically, it is assumed that $Z$ is independent of $Y$ and that, conditionally on $Y=1$ (resp., $Y=-1$), $X$ is a Cox process with intensity $(\lambda_+(t,Z_t))_{t\in [0,T]}$ (resp., $(\lambda_-(t,Z_t))_{t\in [0,T]}$). 

\medskip
It will be assumed that the observation of the trajectories of $X$ is stopped after its $u$-th jump, where $u$ is some known, prespecified, positive integer. Thus, formally, we are to replace $X$ and $Z$ by $X^{\tau}$ and $Z^{\tau}$, where $\tau=\inf\{t \in [0,T]: X_t=u\}$ (stopping time), $X_t^{\tau}=X_{t\wedge \tau}$ and $Z_t^{\tau}=Z_{t\wedge \tau}$. (Notation $t_1 \wedge t_2$ means the minimum of $t_1$ and $t_2$ and, by convention, $\inf \emptyset=0$.) Stopping the observation of $X$ after its $u$-th jump is essentially a technical requirement, with no practical incidence insofar $u$ may be chosen arbitrarily large. However, it should be stressed that with this assumption, $X^{\tau}$ is, with probability one, nicely bounded from above by $u$. Additionally, to keep things simple, we suppose that each $Z_{t}$ takes its values in $[0,1]^d$ and we let $\mathcal Z$ be state space of $Z$.
\medskip

Our objective is to learn the relation between $(X^{\tau},Z^{\tau})$ and $Y$ within the framework
of supervised classification. Given a training dataset of $n$ i.i.d. observation/label pairs $\mathcal D_n=\{(X_1^{\tau_1},Z_1^{\tau_1},Y_1), \hdots, (X_n^{\tau_n},Z_n^{\tau_n},Y_n)\}$ (with evident notation for $\tau_i$'s), distributed as (and independent of) the prototype triplet $(X^{\tau},Z^{\tau},Y)$, the problem is to design a decision rule $g_{n}:\mathcal X\times\mathcal Z\to\{-1,1\}$, based on $\mathcal D_n$, whose role is to assign a label to each possible new instance of the observation $(X^{\tau}, Z^{\tau})$. The classification strategy that we propose is based on empirical convex risk minimization. It is described in the next subsection. 
\medskip

In order to describe our classification procedure,
some more notation is required.
The performance of a classifier $g_{n}:\mathcal X\times\mathcal Z\to\{-1,1\}$ is measured by the probability of error
$$L(g_n)=\mathbb P \left( g_n(X^{\tau},Z^{\tau})\neq Y\,|\,\mathcal D_n\right),$$
and the minimal possible probability of error is the Bayes risk, denoted by
$$L^{\star}=\inf_g L(g)=\mathbb  E \min \left [\eta(X^{\tau}, Z^{\tau}), 1-\eta(X^{\tau}, Z^{\tau})\right].$$
In the identity above, the infimum is taken over all measurable classifiers $g:\mathcal X\times \mathcal Z\to\{-1,  1\}$, and $\eta(X^{\tau},Z^{\tau})=\mathbb P (Y=1\,|\,X^{\tau},Z^{\tau})$
denotes the posterior probability function. The infimum is achieved by the Bayes classifier
$$g^{\star}(X^{\tau},Z^{\tau})=\mbox{sign$\left(2\eta(X^{\tau},Z^{\tau})-1\right)$},$$
where $\mbox{sign$(t)$}=1$ for $t>0$ and $-1$ otherwise.
Our first result (Theorem \ref{th1}) shows that
$$\eta(X^\tau,Z^\tau)=\frac{p_+}{p_-e^{-\xi}+p_+},$$
where $\xi$ is the random variable defined by 
$$\xi=\int_0^{T\wedge \tau}\left(\lambda_--\lambda_+\right)(s,Z_s){\rm d}s+\int_0^{T\wedge \tau} \ln \frac{\lambda_+}{\lambda_-}(s,Z_s){\rm d}X_s.$$

An important consequence is that the Bayes rule associated with our decision problem takes the simple form
$$g^{\star}(X^{\tau},Z^{\tau})=\mbox{sign$\left(\xi-\ln \displaystyle\frac{p_-}{p_+}\right)$}.$$
Next, let $(\varphi_{ j})_{j\ge 1}$ be a countable dictionary of measurable functions defined on $[0,T]\times[0,1]^d$. Assuming that both $\lambda_--\lambda_+$ and $\ln \frac{\lambda_+}{\lambda_-}$ belong to the span of the dictionary, we see that
$$\xi=\sum_{j\ge 1} \left [a_{ j}^{\star}\int_{0}^{T\wedge\tau}\varphi_{ j}(s,Z_{s}){\rm d}s+b^{\star}_{ j}\int_{0}^{T\wedge\tau}\varphi_{ j}(s,Z_{s}){\rm d}X_{s}\right],$$
where $(a^{\star}_j)_{j\geq 1}$ and $(b^{\star}_j)_{j\geq 1}$ are two sequences of unknown real coefficients.
Thus, for each positive integer $B$, it is quite natural to introduce the class $\mathcal F_B$ of real-valued functions $f:\mathcal X\times\mathcal Z\to\R$, defined by
\begin{equation}
\label{FB}
\mathcal F_B=\left\{ f=\sum_{ j=1}^{B} \left [a_{ j} \Phi_{ j}+b_{ j}\Psi_{ j}\right]+c: \max\left(\sum_{j=1}^B|a_{j}|, \sum_{j=1}^B|b_{j}|,|c|\right)\leq B\right\},
\end{equation}
where
$$\Phi_{j}(x,z)=\int_{0}^{T\wedge\tau(x)}\varphi_{ j}(s,z_{s}){\rm d}s,\quad\quad \Psi_{ j}(x,z)=\int_{0}^{T\wedge\tau(x)}\varphi_{ j}(s,z_{s}){\rm d}x_{s},$$
and, by definition, $\tau(x)=\inf\{t\in[0,T]: x_{t}=u\}$ for $x\in\mathcal X$. 
\medskip

Each $f \in \mathcal F_B$ defines a classifier $g_f$ by $g_f=\mbox{sign$(f)$}$.  To simplify notation,
we write $L(f)=L(g_f)=\mathbb P(g_f(X^{\tau},Z^{\tau})\neq Y)$, and note that 
$$ \mathbb E\mathbf 1_{\left[-Yf(X^{\tau},Z^{\tau})>0\right]}\le L(f)\le\mathbb E\mathbf 1_{\left[-Yf(X^{\tau},Z^{\tau})\ge0\right]}.$$
Therefore, the minimization of the probability of error $L(f)$ over $f\in\mathcal F_{B}$ is approximately equivalent to the minimization of the expected 0-1 loss $\mathbf 1_{[.\ge 0]}$ of $-Yf(X^{\tau},Z^{\tau})$. The parameter $B$ may be regarded as an $\mathbb L^1$-type smoothing parameter. Large values of $B$ improve
the approximation properties of the class $\mathcal F_B$ at the price of making the estimation problem
more difficult.
Now, given the sample $\mathcal D_n$, it is reasonable to consider an estimation procedure based on minimizing the sample mean
$$\frac{1}{n}\sum_{i=1}^n\mathbf 1_{[-Y_{i}f(X^{\tau_{i}}_{i},Z^{\tau_{i}}_{i})\ge 0]},$$ 
of the 0-1 loss.
\medskip

It is now well established, however, that such a procedure is computationally intractable as soon as the class $\mathcal F_B$ is nontrivial, since the 0-1 loss function $\mathbf 1_{[.\ge0]}$ is nonconvex. A genuine attempt to circumvent this difficulty is to base the minimization procedure on a convex surrogate $\phi$ of the loss $\mathbf 1_{[.\ge 0]}$. Such convexity-based methods, inspired by the pioneering works on boosting \citep[][]{Fr85,Sc90,FrSc97}, have now largely displaced earlier nonconvex approaches in the machine learning literature \citep[see, e.g.,][and the references therein]{BlLuVa03, LuVa04,Zh04,BaJoAu06}.
\medskip

It turns out that in our Cox process context, the choice of the logit surrogate loss $\phi(t)=\ln_{2}(1+e^t)$ is the most natural one. This will be clarified in Section 2 by connecting the empirical risk minimization procedure and the maximum likelihood principle. Thus, with this choice, the corresponding risk functional and empirical risk functional are defined by
$$A(f)=\mathbb E \phi\left (-Yf(X^{\tau},Z^{\tau})\right)\quad \mbox{and} \quad A_n(f)=\frac{1}{n}\sum_{i=1}^n\phi\left (-Y_if(X_i^{\tau_i},Z_i^{\tau_i})\right).$$
Given a nondecreasing sequence $(B_{k})_{k\ge 1}$ of integer-valued smoothing parameters, the primal estimates we consider take the form
$$\hat f_{k}\in \underset{f \in \mathcal F_{B_{k}}}{\arg\min}\,A_n(f).$$
\begin{rem}
Note that the minimum may not be achieved in $\mathcal F_{B_{k}}$. However, to simplify the arguments, we implicitly assume that the
minimum indeed exists. All proofs may be adjusted, in a straightforward way, to handle approximate minimizers of the empirical risk functional.
\end{rem}
\begin{rem}
The minimization of the functional $A_{n}$ over the class $\mathcal F_{B}$ is indeed a convex problem in the $a_{j}$'s, $b_{j}$'s and $c$ introduced in \eqref{FB}, which makes our method computationally tractable. An alternative approach is to consider functional classes based on the development of the intensity functions $\lambda_{-}$ and $\lambda_{+}$ instead of $\lambda_--\lambda_+$ and $\ln \frac{\lambda_+}{\lambda_-}$. However, such a procedure induces a non convex optimization problem.
\end{rem}

Starting from the collection $(\hat f_k)_{k \geq 1}$, the final estimate uses a value of $k$ chosen empirically, by minimizing a penalized version
of the empirical risk $A_n(\hat f_k)$. To achieve this goal, consider a penalty (or regularization) function $\mbox{pen}:\mathbb N^{\star} \to \mathbb R_{+}$ to be specified later on. Then the resulting penalized estimate $\hat f_{n}=\hat f_{\hat k}$ has
$$\hat k \in \underset{k\ge 1}{\arg\min}\left [ A_n(\hat f_{k})+\mbox{pen$(k)$}\right].$$
The role of the penalty is to compensate for overfitting and helps finding an adequate value of $k$. For larger values of $k$, the class
$\mathcal F_{B_k}$ is larger, and therefore $\mbox{pen$(k)$}$ should be larger as well. 
\medskip

By a careful choice of the regularization term, specified in Theorem \ref{th2},
one may find a close-to-optimal balance between estimation and approximation errors and investigate the probability of error $L(\hat f_n)$ of the classifier $g_{\hat f_n}$ induced by the penalized 
estimate. Our conclusion asserts that $\hat f_n$ adapts nicely to the unknown smoothness of the problem, in the sense that with probability at least $1-1/n^2$,
$$
L(\hat f_{n})-L^{\star}=\mbox{O}\left(\frac{\ln n}{n}\right)^{\frac{\beta}{2\beta+16}},
$$
where $\beta$ is some Sobolev-type regularity measure pertaining to $\lambda_{+}$ and $\lambda_{-}$. 
For the sake of clarity, proofs are postponed to Section 3. An appendix at the end of the paper recalls some important results by \citet{BlBoMa03} and \citet{Ko08} on model selection and suprema of Rademacher processes, together with more technical stochastic calculus material.
\section{Results}
As outlined in the introduction, our first result shows that the posterior probabilities $\mathbb P(Y=\pm 1|X^\tau,Z^\tau)$ have a simple form. The crucial result that is needed here is Lemma \ref{gigi} which uses stochastic calculus arguments. For more clarity, this lemma has been postponed to the Appendix section. Recall that both $p_+$ and $p_-$ are (strictly) positive and satisfy $p_++p_-=1$.
\begin{theo} 
\label{th1}
Let $\xi$ be the  random variable defined by 
$$\xi=\int_0^{T\wedge \tau}(\lambda_--\lambda_+)(s,Z_s){\rm d}s+\int_0^{T\wedge \tau} \ln \frac{\lambda_+}{\lambda_-}(s,Z_s){\rm d}X_s.$$
Then
$$\mathbb P(Y=1|X^\tau,Z^\tau)=\frac{p_+}{p_-e^{-\xi}+p_+} \quad \mbox{and} \quad \mathbb P(Y=-1|X^\tau,Z^\tau)=\frac{p_-}{p_+e^{\xi}+p_-}.$$
\end{theo}

This result, which is interesting by itself, sheds an interesting light on the Cox process classification problem. To see this, fix $Y_1=y_1, \hdots, Y_n=y_n$, and observe that the conditional likelihood of the model is
\begin{align*}
\mathcal L_n & =\prod _{i=1}^n \mathbb P (Y_i=y_i|X^{\tau_i}_i,Z^{\tau_i}_i)\\
& = \prod _{i=1}^n \left(\frac{p_+}{p_-e^{-y_i\xi_i}+p_+}\right)^{\mathbf 1_{[y_i=1]}} \left(\frac{p_-}{p_+e^{-y_i\xi_{i}}+p_-}\right)^{\mathbf 1_{[y_i=-1]}},
\end{align*}
where of course
 $$\xi_i=\int_0^{T\wedge \tau_i}(\lambda_--\lambda_+)(s,Z_{i,s}){\rm d}s+\int_0^{T\wedge \tau_i} \ln \frac{\lambda_+}{\lambda_-}(s,Z_{i,s}){\rm d}X_{i,s}.$$
Therefore, the log-likelihood takes the form
\begin{align*}
\ln \mathcal L_n & = \sum_{i=1}^n \left [\ln \left(\frac{p_+}{p_-e^{-y_i\xi_i}+p_+}\right) {\mathbf 1_{[y_i=1]}} +  \ln \left(\frac{p_-}{p_+e^{-y_i\xi_i}+p_-}\right) {\mathbf 1_{[y_i=-1]}}\right]\\
& =- \sum_{i=1}^n \left [\ln \left(1+\frac{p_-}{p_+}e^{-y_i\xi_i}\right) {\mathbf 1_{[y_i=1]}} +  \ln \left(1+\frac{p_+}{p_-}e^{-y_i\xi_i}\right) {\mathbf 1_{[y_i=-1]}}\right]\\
&=- \sum_{i=1}^n \ln \left(1+\left(\frac{p_-}{p_+}\right)^{y_i}e^{-y_i\xi_i}\right) \\
& =- \sum_{i=1}^n \ln \left(1+\exp \left[ -y_i \left (\xi_i- \ln \frac{p_-}{p_+}\right)\right] \right).
\end{align*}
Thus, letting  $\phi(t)=\ln_{2}(1+e^t)$, we obtain
\begin{equation}
 \label{MLS}
\ln \mathcal L_n = -\ln 2 \sum_{i=1}^n \phi\left(-y_i \left (\xi_i -\ln \frac{p_-}{p_+}\right)\right).
\end{equation}
Since the $\xi_i$'s, $p_+$ and $p_-$ are unknown, the natural idea, already alluded to in the introduction, is to expand $\lambda_--\lambda_+$ and $\ln \frac{\lambda_+}{\lambda_-}$ on the dictionary $(\varphi_j)_{j \geq 1}$. To this end, we introduce the class $\mathcal F_B$ of real-valued functions 
$$\mathcal F_B=\left\{ f=\sum_{ j=1}^{B} \left [a_{ j} \Phi_{ j}+b_{ j}\Psi_{ j}\right]+c: \max\left(\sum_{j=1}^B|a_{j}|, \sum_{j=1}^B|b_{j}|,|c|\right)\leq B\right\},$$
where $B$ is a positive integer, 
$$\Phi_{j}(x,z)=\int_{0}^{T\wedge\tau(x)}\varphi_{ j}(s,z_{s}){\rm d}s,\quad\mbox{and}\quad \Psi_{ j}(x,z)=\int_{0}^{T\wedge\tau(x)}\varphi_{ j}(s,z_{s}){\rm d}x_{s}.$$
For a nondecreasing sequence $(B_{k})_{k\ge1}$ of integer-valued smoothing parameters and for each $k\ge 1$, we finally select $\hat f_k\in \mathcal F_{B_{k}}$ for which the log-likelihood (\ref{MLS}) is maximal. Clearly, such a maximization strategy is strictly equivalent to minimizing over $f\in \mathcal F_{B_{k}}$ the empirical risk 
$$A_n(f)=\frac{1}{n}\sum_{i=1}^n\phi\left (-Y_if(X_i^{\tau_i},Z_i^{\tau_i})\right).$$
This remark reveals the deep connection between our Cox process learning model and the maximum likelihood principle. In turn, it justifies the logit loss $\phi(t)=\ln_{2}(1+e^t)$ as the natural surrogate candidate to the nonconvex 0-1 classification loss. (Note that the $\ln 2$ term is introduced for technical reasons only and plays no role in the analysis). Finally, we stress the fact that by convexity, this approach is computationally tractable. 
\medskip

As for now, denoting by $\|.\|_{\infty}$ the functional supremum norm, we assume that there exists a positive constant $L$ such that, for each $j\geq 1$, $\|\varphi_{j}\|_{\infty}\leq L$. It immediately follows that for all integers $B\ge 1$, the class $\mathcal F_{B}$ is uniformly bounded by $UB$, where $U=1+(T+u)L$. We are now ready to state our main theorem, which offers a bound on the difference $A(\hat f_{n})-A(f^{\star})$.
\begin{theo}
\label{th2}
Let $(B_{k})_{k\geq 1}$ be a nondecreasing sequence of positive integers such that $\sum_{k\geq 1}B^{-\alpha}_{k}\leq 1$ for some $\alpha>0$. For all $k\ge 1$, let
$$R_{k}=A^2_{k}B_{k}C_{k}+\frac{\sqrt{A_{k}}}{C_{k}},$$ 
where 
$$A_{k}=UB_{k}\phi'(UB_{k}) \quad \mbox{and} \quad C_{k}=2(\phi(UB_{k})+1-\ln 2).$$
Then there exists a universal constant $C>0$ such that if the penalty ${\rm pen}:\N^{\star}\ra\R_{+}$ satisfies
$$\emph{pen$(k)$}\ge C\left[R_{k}\frac{\ln n}{n}+\frac{C_{k}(\alpha\ln B_{k}+\delta+\ln 2)}{n}\right]$$
for some $\delta> 0$, one has, with probability at least $1-e^{-\delta}$,
\begin{equation}
\label{IO}
A(\hat f_{n})-A(f^{\star})\leq 2\inf_{k\ge 1}\left\{\inf_{f\in\mathcal F_{B_{k}}}\left(A(f)-A(f^{\star})\right)+\emph{pen$(k)$}\right\}.
\end{equation}
\end{theo}

Some remarks are in order. At first, we note that Theorem \ref{th2} provides us with an oracle inequality which shows that, for each $B_k$, the penalized estimate does almost as well as the best possible classifier in the class $\mathcal F_{B_k}$, up to a term of the order $\ln n/n$. It is stressed that this remainder term tends to 0 at a much faster rate than the standard $(1/\sqrt n)$-term suggested by a standard uniform convergence argument \citep[see, e.g.,][]{LuVa04}. This is a regularization effect which is due to the convex loss $\phi$. In fact, proof of Theorem \ref{th2} relies on the powerful model selection machinery presented in \citet{BlBoMa03} coupled with modern empirical process theory arguments developed in \citet{Ko08}. We also emphasize that a concrete but suboptimal value of the constant $C$ may be deduced from the proof, but that no attempt has been made to optimize this constant. 
Next, observing that, for the logit loss, 
$$\phi'(t)=\frac{1}{\ln 2(e^{-t}+1)},$$
we notice that a penalty behaving as $B^{4}_k$ is sufficient for the oracle inequality of Theorem \ref{th2} to hold. This corresponds to a regularization function proportional to the fourth power of the $\mathbb L^1$-norm of the collection of coefficients defining the base class functions. Such regularizations have been explored by a number of authors in recent years, specifically in the context of sparsity and variable selection \citep[see, e.g.,][]{Ti96,Ca05, BuTsWe07,BiRi09}. With this respect, our approach is close to the view of \citet{MaMe11}, who provide information about the Lasso as an $\mathbb L^1$-regularization procedure per se, together with sharp $\mathbb L^1$-oracle inequalities. Let us finally mention that the result of Theorem \ref{th2} can be generalized, with more technicalities, to other convex loss functions by following, for example, the arguments presented in \citet{BaJoAu06}. 
\medskip

If we are able to control the approximation term $\inf_{f\in\mathcal F_{B_{k}}}(A(f)-A(f^{\star}))$ in inequality (\ref{IO}), then it is possible to give an explicit rate of convergence to 0 for the quantity $A(\hat f_{n})-A(f^{\star})$. This can be easily achieved by assuming, for example, that $(\varphi_j)_{j \geq 1}$ is an orthonormal basis and that both combinations $\lambda_{-}-\lambda_{+}$ and $\ln\frac{\lambda_{+}}{\lambda_{-}}$ enjoy some Sobolev-type regularity with respect to this basis. Also, the following additional assumption will be needed:
\medskip

\noindent\textbf{Assumption $\mathbf A$.} There exists a measure $\mu$ on $[0,1]^d$ and a constant $D>0$ such that, for all $t\in[0,T]$, the distribution of $Z_{t}$ has a density $h_{t}$ with respect to $\mu$ which is uniformly bounded by $D$. In addition, $\lambda_{-}$ and $\lambda_{+}$ are both $[\e,D]$-valued for some $\e>0$. 
\begin{pro}
\label{pro4}
Assume that Assumption $\mathbf A$ holds. Assume, in addition, that $(\varphi_{j})_{j\ge 1}$ is an orthonormal basis of $\mathbb L^2({\rm d}s\otimes\mu)$, where ${\rm d}s$ stands for the Lebesgue measure on $[0,T]$, and that both $\lambda_{-}-\lambda_{+}$ and $\ln\frac{\lambda_{+}}{\lambda_{-}}$ belong to the ellipso\"{i}d
$$\mathcal W(\beta,M)=\left\{f=\sum_{j=1}^{\infty}a_{j}\varphi_{j}:\sum_{j=1}^{\infty}j^{2\beta}a^{2}_{j}\leq M^2\right\},$$
for some fixed $\beta\in\mathbb N^{\star}$ and $M>0$. Then, letting
$$\lambda_{-}-\lambda_{+}=\sum_{j=1}^{\infty}a^{\star}_{j}\varphi_{j}\quad\mbox{and}\quad\ln\frac{\lambda_{+}}{\lambda_{-}}=\sum_{j=1}^{\infty}b^{\star}_{j}\varphi_{j},$$
we have, for all $B\ge\max(M^2, \ln\frac{p_{+}}{p_{-}})$,
\begin{align*}
\inf_{f\in\mathcal F_{B}}\left(A(f)-A(f^{\star})\right) &\le  \frac{2D\sqrt{T\mu\left([0,1]^d\right)M\| a^{\star}\|_{2}}}{B^{\beta/2}}
\nonumber\\
& \quad +\ \frac{2D(1+D\sqrt{T\mu([0,1]^{d})})\sqrt{M\| b^{\star}\|_{2}}}{B^{\beta/2}},
\nonumber
\end{align*}
where $\| a^{\star}\|^2_{2}=\sum_{j= 1}^{\infty} a^{\star2}_{j}$ and $\|b^{\star}\|^{2}_{2}=\sum_{j=1}^{\infty}b^{\star2}_{j}$.
\end{pro}
A careful inspection of Theorem \ref{th2} and Proposition \ref{pro4} reveals that for the choice $B_{k}=\lceil (\pi k)^{2/\alpha}/6^{1/\alpha}\rceil$ and $\delta=2\ln n$, there exists a universal constant $C>0$ such that
$$A(\hat f_{n})-A(f^{\star})\leq C\mathbf L\left(\sqrt{\Vert a^{\star}\Vert_{2}}+\sqrt{\Vert b^{\star}\Vert_{2}}\right)^{\frac{8}{\beta+8}}\left(\frac{\ln n}{n}\right)^{\frac{\beta}{\beta+8}},$$
with probability at least $1-1/n^2$, where 
$$\mathbf L=U^{\frac{3\beta}{\beta+8}}\left[2D\sqrt{M}(1+D\sqrt{T\mu([0,1]^{d})})\right]^{\frac{8}{\beta+8}}\max\left(\left(\frac{\beta}{8}\right)^{\frac{8}{\beta+8}},\left(\frac{8}{\beta}\right)^{\frac{\beta}{\beta+8}}\right).$$
Observe that, due to the specific form of the ellipso\"{i}d $\mathcal W(\beta,M)$, the rate of convergence does not depend upon the dimension $d$.
\medskip

Of course, our main concern is not the behavior of the expected risk $A(\hat f_n)$ but the probability of error $L(\hat f_n)$ of the corresponding classifier. Fortunately, the difference $L(\hat f_n)-L^{\star}$ may directly be related to $A(\hat f_n)-A(f^{\star})$. Applying for example Lemma 2.1 in \citet{Zh04}, we conclude that with probability at least $1-1/n^2$,
\begin{equation}
\label{AB}
L(\hat f_{n})-L^{\star}\leq 2\sqrt{C\mathbf L}\left(\sqrt{\Vert a^{\star}\Vert_{2}}+\sqrt{\Vert b^{\star}\Vert_{2}}\right)^{\frac{4}{\beta+8}}\left(\frac{\ln n}{n}\right)^{\frac{\beta}{2\beta+16}}.
\nonumber
\end{equation}
To understand the significance of this inequality, just recall that what we are after in this article is the supervised classification of (infinite-dimensional) stochastic processes. As enlightened in the proofs, this makes the analysis different from the standard context, where one seeks to learn finite-dimensional quantities. The bridge between the two worlds is crossed via  stochastic calculus arguments. Lastly, it should be noted that the regularity parameter $\beta$ is assumed to be unknown, so that our results are adaptive as well.

\section{Proofs}
Throughout this section, if $P$ is a probability measure and $f$ a function, the notation $Pf$ stands for the integral of $f$ with respect to $P$. By $\mathbb L^2(P)$ we mean the space of square integrable real functions with respect to $P$. Also, for a class $\mathcal F$ of functions in $\mathbb L^2(P)$ and $\e>0$, we denote by $N(\e,\mathcal F,\mathbb L^2(P))$ the $\e$-covering number of $\mathcal F$ in $\mathbb L^2(P)$, i.e., the minimal number of metric balls of radius $\e$ in $\mathbb L^2(P)$ that are needed to cover $\mathcal F$ \citep[see, e.g., Definition 2.1.5 in][]{vaWe96}.

\subsection{Proof of Theorem \ref{th1}}
For any stochastic processes $M_1$ and $M_2$, the notation $\mathbb Q_{M_2|M_1}$ and $\mathbb Q_{M_2}$ respectively mean the distribution under $\mathbb Q$ of $M_2$ given $M_1$, and the distribution under $\mathbb Q$ of $M_2$. 
\medskip

We start the proof by observing that
\begin{equation}\label{densiteRN}
\mathbb P(Y=1\,|\,X^\tau=x,Z=z)=p_+ \frac{{\rm d} \mathbb P_{X^\tau,Z|Y=1}}{{\rm d} \mathbb P_{X^\tau,Z}}(x,z).
\end{equation}
Thus, to prove the theorem, we need to evaluate the above Radon-Nikodym density. To this aim, we introduce the conditional probabilities $\mathbb P^{\pm}=\mathbb P(.|Y=\pm 1)$. For any path $z$ of $Z$, the conditional distributions $\mathbb P^+_{X|Z=z}$ and $\mathbb P^-_{X|Z=z}$  are those of Poisson processes with intensity $\lambda_+(.,z)$ and $\lambda_-(.,z)$, respectively. Consequently, according to Lemma \ref{gigi}, the stopped process $X^\tau$ satisfies 
$$D_+(x,z)\,  \mathbb P^+_{X^\tau|Z=z}({\rm d}x)=D_-(x,z)\, \mathbb P^-_{X^\tau|Z=z}({\rm d}x),$$
where 
$$D_\pm(x,z)=\exp \left(-\int_0^{T\wedge \tau}\left(1-\lambda_\pm(s,z_s)\right){\rm d}s-\int_0^{T\wedge \tau}\ln \lambda_\pm (s,z_s){\rm d}x_s\right).$$
Therefore, 
$$D_+(x,z)\, \mathbb P^+_{X^\tau|Z=z}\otimes \mathbb P_Z ({\rm d}x,{\rm d}z)=D_-(x,z)\, \mathbb P^-_{X^\tau|Z=z}\otimes \mathbb P_Z ({\rm d}x,{\rm d}z).$$
But, by independence of $Y$ and $Z$, one has $\mathbb P_Z=\mathbb P^+_Z=\mathbb P^-_Z$. Thus, 
$$\mathbb P_{X^\tau,Z|Y=\pm 1}({\rm d}x,{\rm d}z)=\mathbb P^\pm_{X^\tau|Z=z}\otimes \mathbb P_Z ({\rm d}x,{\rm d}z),$$
whence
$$D_+(x,z)\mathbb P_{X^\tau,Z|Y=1}({\rm d}x,{\rm d}z)=D_-(x,z)\mathbb P_{X^\tau,Z|Y=-1}({\rm d}x,{\rm d}z).$$
On the other hand,
$$\mathbb P_{X^\tau,Z}(x,z)=p_+ \mathbb P_{X^\tau,Z|Y=1}(x,z)+p_- \mathbb P_{X^\tau,Z|Y=-1}(x,z),$$
so that
$$\frac{{\rm d} \mathbb P_{X^\tau,Z|Y=1}}{{\rm d} \mathbb P_{X^\tau,Z}}(x,z)=\frac{1}{p_- \frac{D_+(x,z)}{D_-(x,z)}+p_+}.$$
Using identity \eqref{densiteRN}, we obtain
$$\mathbb P(Y=1\,|\,X^\tau,Z)=\frac{p_+}{p_-e^{-\xi}+p_+},$$
where
$$
\xi  = \int_0^{T\wedge \tau}(\lambda_--\lambda_+)(s,Z_s){\rm d}s+\int_0^{T\wedge \tau} \ln \frac{\lambda_+}{\lambda_-}(s,Z_s){\rm d}X_s.$$
Observing now that $\sigma(\tau)\subset \sigma(X^\tau_t, t\le T)$ and 
$$\xi =  \int_0^{T\wedge \tau}(\lambda_--\lambda_+)(s,Z_s^\tau){\rm d}s+\int_0^{T\wedge \tau} \ln \frac{\lambda_+}{\lambda_-}(s,Z_s^\tau){\rm d}X_s^\tau$$
give 
$$\mathbb P(Y=1\,|\,X^\tau,Z^\tau)=\mathbb P(Y=1\,|\,X^\tau,Z).$$
This shows the desired result.$\Box$

\subsection{Proof of Theorem \ref{th2}}

Theorem \ref{th2} is mainly a consequence of a general model selection result due to \citet{BlBoMa03}, which is recalled in the Appendix for the sake of completeness (Theorem \ref{B}). Throughout the proof, the letter $C$ denotes a generic universal positive constant, whose value may change from line to line. We let $\ell(f)$ be a shorthand notation for the function 
$$(x,z,y)\in\mathcal X\times\mathcal Z\times\{-1,1\}\mapsto\phi(-yf(x,z)),$$
and let $P$ be the distribution of the prototype triplet $(X^{\tau},Z^{\tau},Y)$.
\medskip

To frame our problem in the vocabulary of Theorem \ref{B}, we consider the family of models $(\mathcal F_{B_{k}})_{k\geq 1}$ and start by verifying that assumptions $(i)$ to $(iv)$ are satisfied. If we define
$$\mathbf d^2(f,f')=P\left(\ell(f)-\ell(f')\right)^2,$$
then assumption $(i)$ is immediately satisfied. A minor modification of the proof of Lemma 19 in \citet{BlLuVa03} reveals that, for all integers $B>0$ and all $f\in\mathcal F_{B}$, 
$$P\left(\ell(f)-\ell(f^{\star})\right)^2\leq \left(\phi(UB)+\phi(-UB)+2-2\ln2\right)P\left(\ell(f)-\ell(f^{\star})\right).$$ 
This shows that assumption $(ii)$ is satisfied with $C_{k}=2(\phi(UB_{k})+1-\ln 2)$. Moreover, it can be easily verified that assumption $(iii)$ holds with $b_{k}=\phi(UB_{k})$. 
\medskip

The rest of the proof is devoted to the verification of assumption $(iv)$. To this aim, for all $B>0$ and all $f_{0}\in\mathcal F_{B}$, we need to bound the expression
$$F_{B}(r)=\esp\sup\left\{\left\vert(P_{n}-P)\left(\ell(f)-\ell(f_{0})\right)\right\vert: f\in\mathcal F_{B}, \mathbf d^2(f,f_{0})\leq r\right\},$$
where
$$P_{n}=\frac{1}{n}\sum_{i=1}^{n}\delta_{\left(X^{\tau_{i}}_{i},Z^{\tau_{i}}_{i},Y_{i}\right)}$$ 
is the empirical distribution associated to the sample. Let 
$$\mathcal G_{B,f_{0}}=\big\{\ell(f)-\ell(f_{0}):f\in\mathcal F_{B}\big\}.$$
Then
$$F_{B}(r)=\esp\sup\left\{\left |(P_{n}-P)g\right |: g\in\mathcal G_{B,f_{0}}, Pg^2\leq r\right\}.$$
Using the symmetrization inequality presented in Theorem 2.1 of \citet{Ko08}, it is easy to see that
\begin{equation}
\label{e1}
F_{B}(r)\leq 2\,\esp\sup\left\{\frac{1}{n}\sum_{i=1}^{n}\sigma_{i}g\left(X^{\tau_{i}}_{i},Z^{\tau_{i}}_{i},Y_{i}\right): g\in\mathcal G_{B,f_{0}}, Pg^2\leq r\right\},
\end{equation}
where $\sigma_{1},\dots,\sigma_{n}$ are independent Rademacher random variables (that is, $\prob(\sigma_{i}=\pm1)=1/2$), independent from the $(X^{\tau_{i}}_{i},Z^{\tau_{i}}_{i},Y_{i})$'s. Now, since the functions in $\mathcal F_{B}$ take their values in $[-UB,UB]$, and since $\phi$ is Lipschitz on this interval with constant $\phi'(UB)$, we have, for all $f,f'\in\mathcal F_{B}$, 
$$\sqrt{P_{n}\left(\ell(f)-\ell(f')\right)^2}\le\phi'(UB)\sqrt{P_{n}\left(f-f'\right)^2}.$$
Consequently, for all $\e>0$, 
$$N\left(2\e UB\phi'(UB),\mathcal G_{B,f_{0}},\mathbb L^2\left(P_{n}\right)\right)\leq N\left(2\e UB,\mathcal F_{B},\mathbb L^2\left(P_{n}\right)\right).$$
Since $\mathcal F_{B}$ is included in a linear space of dimension at most $2B+1$, Lemma 2.6.15 in \citet{vaWe96} indicates that it is a VC-subgraph class of VC-dimension at most $2B+3$. Observing that the function constantly equal to $2UB$ is a measurable envelope for $\mathcal F_{B}$, we conclude from Theorem 9.3 in \citet{Kos08} that, for all $\e>0$, 
$$N\left(2\e UB,\mathcal F_{B},\mathbb L^2\left(P_{n}\right)\right)\leq C\left(2B+3\right)\left(4e\right)^{2B+3}\left(\frac{1}{\e}\right)^{4\left(B+1\right)}.$$
Therefore, 
$$N\left(2\e UB\phi'(UB),\mathcal G_{B,f_{0}},\mathbb L^2\left(P_{n}\right)\right)\leq C\left(2B+3\right)\left(4e\right)^{2B+3}\left(\frac{1}{\e}\right)^{4\left(B+1\right)}.$$
Now, notice that the constant function equal to $2UB\phi'(UB)$ is a measurable envelope for $\mathcal G_{B,f_{0}}$. Thus, applying Lemma \ref{A} yields 
$$F_{B}(r)\le\psi_{B}(r),$$
where $\psi_{B}$ is defined for all $r>0$ by
$$\psi_{B}(r)=\frac{C\sqrt r}{\sqrt n}\sqrt{B\ln\left(\frac{A'_{B}}{\sqrt r}\right)}\lor\frac{CBA_{B}}{ n}\ln\left(\frac{A'_{B}}{\sqrt r}\right)\lor\frac{CA_{B}}{n}\sqrt{B\ln\left(\frac{A'_{B}}{\sqrt r}\right)},$$
with $A_{B}=UB\phi'(UB)$ and $A'_{B}=A_{B}((2B+3)(4e)^{2B+3})^{1/4(B+1)}$. (Notation $t_1 \lor t_2$ means the maximum of $t_1$ and $t_2$.) 
\medskip

Attention shows that $\psi_{B}$ is a sub-root function and assumption $(iv)$ is therefore satisfied. It is routine to verify that the solution $r^{\star}_{k}$ of $\psi_{B_{k}}(r)=r/C_{k}$ satisfies, for all $k\ge 1$ and all $n\ge 1$,
$$r^{\star}_{k}\leq C\left(A^{2}_{B_{k}}B_{k}C^2_{k}+\sqrt {A'_{B_{k}}}\right)\frac{\ln n}{n}.$$
Furthermore, observing that the function $B\mapsto((2B+3)(4e)^{2B+3})^{1/4(B+1)}$ is bounded from above, we obtain
$$r^{\star}_{k}\leq C\left(A^{2}_{B_{k}}B_{k}C^2_{k}+\sqrt {A_{B_{k}}}\right)\frac{\ln n}{n}.$$
Hence, taking $x_{k}=\alpha\ln \lambda_{k}$ and $K=11/5$ in Theorem \ref{B}, and letting
$$R_{k}=A^2_{B_{k}}B_{k}C_{k}+\frac{\sqrt{A_{B_{k}}}}{C_{k}},$$
we conclude that there exists a universal constant $C>0$ such that, if the penalty $\mbox{pen}:\N^{\star}\ra\R_{+}$ satisfies
$$\mbox{pen$(k)$}\ge C\left\{R_{k}\frac{\ln n}{n}+\frac{C_{k}\left(\alpha\ln B_{k}+\delta+\ln 2\right)}{n}\right\}$$
for some $\delta>0$, then, with probability at least $1-e^{-\delta}$, 
$$A(\hat f_{n})-A(f^{\star})\leq 2\inf_{k\ge 1}\left\{\inf_{f\in\mathcal F_{B_{k}}}(A(f)-A(f^{\star}))+\mbox{pen$(k)$}\right\}.$$
This completes the proof. $\square$
\subsection{Proof of Proposition \ref{pro4}}
Proof of Proposition \ref{pro4} relies on the following intermediary lemma,
which is proved in the next subsection.
\begin{lem}
\label{lem3}
Assume that Assumption $\mathbf A$ holds. Then, for all positive integers $B\ge 1$, 
\begin{align*}
 \inf_{f\in\mathcal F_{B}}\left(A(f)-A(f^{\star})\right) & \le  2D \sqrt{T\mu([0,1]^{d})}\min\left\|\sum_{j=1}^{B}\alpha_{j}\varphi_{j}-\left(\lambda_{-}-\lambda_{+}\right)\right\|\\
& \quad + 2D(1+D\sqrt{T\mu([0,1]^{d})})\min \left\|\sum_{j=1}^{B}\alpha_{j}\varphi_{j}-\ln\frac{\lambda_{+}}{\lambda_{-}}\right\|\\
& \quad + 2\min_{\vert x\vert\leq B}\left \vert x-\ln\frac{p_{+}}{p_{-}}\right\vert,
\end{align*}
where the first two minima are taken over all $\alpha=(\alpha_{1},\dots,\alpha_{B})\in\R^B$ with $\sum_{j=1}^B|\alpha_j|\leq B$ and where we have denoted by $\Vert .\Vert$ the $\mathbb L^{2}({\rm d}s\otimes\mu)$-norm. 
\end{lem}

\noindent{\sc Proof of Proposition \ref{pro4}} -- For ease of notation, we will denote by $\Vert .\Vert$ the $\mathbb L^{2}({\rm d}s\otimes\mu)$-norm throughout the proof. For all $B\ge 1$,
\begin{align}
&\min\left\{\left\|\sum_{j=1}^{B}\alpha_{j}\varphi_{j}-\left(\lambda_{-}-\lambda_{+}\right)\right\|:\sum_{j=1}^{B}| \alpha_{j}|\leq B\right\}
\nonumber\\
&\quad\leq \min\left\{\left\|\sum_{j=1}^{B}\alpha_{j}\varphi_{j}-\left(\lambda_{-}-\lambda_{+}\right)\right\|:\sum_{j=1}^{B}\alpha^2_{j}\leq B\right\}.
\label{ref-cor4-2}
\end{align}
Since $\lambda_{-}-\lambda_{+}\in\mathcal W(\beta,M)$ and $B\ge M^2$, we have 
\begin{equation}
\label{ref-cor4-3}
\sum_{j=1}^{B}a^{\star2}_{j}\le\sum_{j=1}^{\infty}j^{2\beta}a^{\star2}_{j}\leq M^2\leq B.
\end{equation}
Thus, combining \eqref{ref-cor4-2} and \eqref{ref-cor4-3} yields, for $B\ge M^2$,
\begin{align}
\min\left\{\left\|\sum_{j=1}^{B}\alpha_{j}\varphi_{j}-\left(\lambda_{-}-\lambda_{+}\right)\right\|:\sum_{j=1}^{B}| \alpha_{j}|\leq B\right\}
&\le \left\|\sum_{j=1}^{B}a^{\star}_{j}\varphi_{j}-\left(\lambda_{-}-\lambda_{+}\right)\right\|
\nonumber\\
&=\left\|\sum_{j=B+1}^{\infty}a^{\star}_{j}\varphi_{j}\right\|.
\label{ref-cor4-4}
\end{align}
It follows from the properties of an orthonormal basis and the definition of $\mathcal W(\beta, M)$ that
\begin{align}
\left\|\sum_{j=B+1}^{\infty}a^{\star}_{j}\varphi_{j}\right\|^2& =\sum_{j=B+1}^{\infty}a^{\star2}_{j}
\nonumber\\
&\leq \sqrt{\sum_{j=B+1}^{\infty}j^{2\beta}a^{\star2}_{j}}\sqrt{\sum_{j=B+1}^{\infty}\frac{a^{\star2}_{j}}{j^{2\beta}}}
\nonumber\\
&\leq M\sqrt{\sum_{j=B+1}^{\infty}\frac{a^{\star2}_{j}}{j^{2\beta}}}
\nonumber\\
&\leq \frac{M\| a^{\star}\|_{2}}{B^\beta}.
\label{ref-cor4-5}
\end{align}
Inequalities \eqref{ref-cor4-4} and \eqref{ref-cor4-5} show that, for all $B\ge M^2$, 
$$
\min\left\{\left\|\sum_{j=1}^{B}\alpha_{j}\varphi_{j}-\left(\lambda_{-}-\lambda_{+}\right)\right\|:\sum_{j=1}^{B}| \alpha_{j}|\leq B\right\}\leq \sqrt{\frac{M\| a^{\star}\|_{2}}{B^\beta}}.
$$

Similarly, it may be proved that, for all $B\ge M^2$, 
$$
\min\left\{\left\|\sum_{j=1}^{B}\alpha_{j}\varphi_{j}-\ln\frac{\lambda_{+}}{\lambda_{-}}\right\|:\sum_{j=1}^{B}| \alpha_{j}|\leq B\right\}\leq \sqrt{\frac{M\| b^{\star}\|_{2}}{B^\beta}}.
$$
Applying Lemma \ref{lem3} we conclude that, whenever $B\ge\max(M^2, \ln\frac{p_{+}}{p_{-}})$, 
\begin{align*}
\inf_{f\in\mathcal F_{B}}\left(A(f)-A(f^{\star})\right) &\le  \frac{2D\sqrt{T\mu\left([0,1]^d\right)M\| a^{\star}\|_{2}}}{B^{\beta/2}}
\nonumber\\
& \quad+ \frac{2D(1+D\sqrt{T\mu([0,1]^{d})})\sqrt{M\| b^{\star}\|_{2}}}{B^{\beta/2}},
\nonumber
\end{align*}
which ends the proof. $\square$
\subsection{Proof of Lemma \ref{lem3}}
We start with a technical lemma.
\begin{lem}
\label{lem1}
Let $\phi(t)=\ln_{2}(1+e^t)$ be the logit loss. Then
$$\underset{f}{\arg\min}\,\esp\phi\left(-Yf(X^\tau,Z^\tau)\,|\, X^\tau,Z^\tau\right)=\xi-\ln\frac{p_{-}}{p_{+}},$$
where the minimum is taken over all measurable functions $f:\mathcal X\times\mathcal Z\ra\R.$
\end{lem}

\noindent{\sc Proof} -- According to the results of Section $2.2$ in \citet{BaJoAu06}, one has
$$\underset{f}{\arg\min}\ \esp\phi\left(-Yf(X^\tau,Z^\tau)\,|\,X^\tau,Z^\tau\right)=\alpha^{\star}\left(\eta(X^\tau,Z^\tau)\right),$$
where, for all $0\le\eta\leq 1$, 
$$\alpha^{\star}(\eta)=\underset{\alpha\in\R}{\arg\min}\left(\eta\phi(-\alpha)+(1-\eta)\phi(\alpha)\right).$$
With our choice for $\phi$, it is straightforward to check that, for all $0\le\eta<1$, 
$$\alpha^{\star}(\eta)=\ln\left(\frac{\eta}{1-\eta}\right).$$
Since, by assumption, $p_{-}>0$, we have
$$\eta(X^\tau,Z^\tau)=\frac{p_+}{p_-e^{-\xi}+p_+}<1.$$
Thus  
$$\alpha^{\star}\left(\eta(X^\tau,Z^\tau)\right)=\xi-\ln\frac{p_{-}}{p_{+}},$$
which is the desired result. $\square$
\medskip

\noindent{\sc Proof of Lemma \ref{lem3}} -- Let $B>0$ be fixed. Let $a_{1},\dots,a_{B}$ and $b_{1},\dots,b_{B}$ be real numbers such that
$$\left \|\sum_{j=1}^{B}a_{j}\varphi_{j}-\left(\lambda_{-}-\lambda_{+}\right)\right\|_{\mathbb L^{2}({\rm d}s\otimes\mu)}=\min\left\|\sum_{j=1}^{B}\alpha_{j}\varphi_{j}-\left(\lambda_{-}-\lambda_{+}\right)\right\|_{\mathbb L^{2}({\rm d}s\otimes\mu)}$$
and
$$\left\|\sum_{j=1}^{B}b_{j}\varphi_{j}-\ln\frac{\lambda_{+}}{\lambda_{-}}\right\|_{\mathbb L^{2}({\rm d}s\otimes\mu)}=\min \left\|\sum_{j=1}^{B}\alpha_{j}\varphi_{j}-\ln\frac{\lambda_{+}}{\lambda_{-}}\right\|_{\mathbb L^{2}({\rm d}s\otimes\mu)},$$
where, in each case, the minimum is taken over all $\alpha=(\alpha_{1},\hdots,\alpha_{B})\in\R^{B}$ with $\sum_{j=1}^B| \alpha_{j}|\leq B$. Let also $c\in\R$ be such that
$$\left\vert c-\ln\frac{p_{+}}{p_{-}}\right\vert=\min_{\vert x\vert\leq B}\left\vert x-\ln\frac{p_{+}}{p_{-}}\right\vert.$$
Introduce $f_{B}$, the function in $\mathcal F_{B}$ defined by 
\begin{align}
f_{B} &=\sum_{j=1}^B\left[ a_{j}\Phi_{j}+b_{j}\Psi_{j}\right]+c
\nonumber\\
& =\int_{0}^{T\wedge\tau}\sum_{j=1}^B a_{j}\varphi_{j}(s,Z_{s}){\rm d}s+\int_{0}^{T\wedge\tau}\sum_{j=1}^B b_{j}\varphi_{j}(s,Z_{s}){\rm d}X_{s}+c.
\nonumber
\end{align}
Clearly,
\begin{equation}
\label{ref-th3-0}
\inf_{f\in\mathcal F_{B}}\left(A(f)-A(f^{\star})\right)\leq A(f_{B})-A(f^{\star}).
\end{equation}
Since $\phi$ is Lipschitz with constant $\phi'(UB)=(\ln2(1+e^{-UB}))^{-1}\leq 2$ on the interval $[-UB,UB]$, we have
\begin{equation}
\label{ref-th3-1}
\left| A(f_{B})-A(f^{\star})\right|\leq 2\esp \left | f_{B}(X^\tau,Z^\tau)-f^{\star}(X^\tau,Z^\tau)\right|.
\end{equation}
But, by Lemma \ref{lem1}, 
$$f^{\star}(X^{\tau},Z^{\tau})=\int_0^{T\wedge \tau}\left(\lambda_--\lambda_+\right)(s,Z_s){\rm d}s+\int_0^{T\wedge \tau} \ln \frac{\lambda_+}{\lambda_-}(s,Z_s){\rm d}X_s+\ln\frac{p_{+}}{p_{-}}.$$
Thus, letting, 
$$\vartheta_{1}=\sum_{j=1}^B a_{j}\varphi_{j}-\left(\lambda_--\lambda_+\right)\quad\mbox{and}\quad \vartheta_{2}=\sum_{j=1}^B b_{j}\varphi_{j}-\ln \frac{\lambda_+}{\lambda_-},$$
it follows 
\begin{align}
\esp \left | f_{B}(X^\tau,Z^\tau)-f^{\star}(X^\tau,Z^\tau)\right|\ & \leq \esp\left| \int_0^{T\wedge \tau} \vartheta_{1}(s,Z_{s}){\rm d}s \right|
\nonumber\\
&\quad+ \esp\left| \int_0^{T\wedge \tau} \vartheta_{2}(s,Z_{s}){\rm d}X_{s} \right|
\nonumber\\
&\quad+ \left| c-\ln\frac{p_{+}}{p_{-}}\right|.
\label{ref-th3-2}
\end{align}
Using Assumption $\mathbf A$ and Cauchy-Schwarz's Inequality, we obtain
\begin{align}
\esp\left | \int_0^{T\wedge \tau} \vartheta_{1}(s,Z_{s}){\rm d}s \right | & \leq \int_0^{T}\int_{[0,1]^{d}} \left |\vartheta_{1}(s,z)\right |\prob_{Z_{s}}({\rm d}z){\rm d}s 
\nonumber\\
& =\int_0^{T}\int_{[0,1]^{d}} \left |\vartheta_{1}(s,z)\right | h_{s}(z)\mu({\rm d}z){\rm d}s 
\nonumber\\
& \leq D\| \vartheta_{1}\|_{\mathbb L^{1}({\rm d}s\otimes\mu)}
\nonumber\\
& \leq D\sqrt{T\mu([0,1]^{d})}\| \vartheta_{1}\|_{\mathbb L^{2}({\rm d}s\otimes\mu)}.
\label{ref-th3-3}
\end{align}
With a slight abuse of notation, set $\lambda_{Y}=\lambda_{\pm}$, depending on whether $Y=\pm1$, and
$$\Lambda_{Y,Z}(t)=\int_{0}^{t}\lambda_{Y}(s,Z_{s}){\rm d}s, \quad t\in[0,T].$$ 
With this notation,
\begin{align}
\esp\left | \int_0^{T\wedge \tau} \vartheta_{2}(s,Z_{s}){\rm d}X_{s} \right | & \leq \esp\left | \int_0^{T\wedge\tau} \vartheta_{2}(s,Z_{s}){\rm d}\left(X_{s}-\Lambda_{Y,Z}(s)\right) \right | 
\nonumber\\
&\quad + \esp\left | \int_0^{T\wedge\tau} \vartheta_{2}(s,Z_{s}){\rm d}\Lambda_{Y,Z}(s) \right | 
\nonumber\\
&=\esp\left | \int_0^{T\wedge\tau} \vartheta_{2}(s,Z_{s}){\rm d}\left(X_{s}-\Lambda_{Y,Z}(s)\right) \right | 
\nonumber\\
& \quad+ \esp\left | \int_0^{T\wedge\tau} \vartheta_{2}(s,Z_{s})\lambda_{Y}(s,Z_{s}){\rm d}s\right |. 
\label{ref-th3-4}
\end{align}
Therefore, applying Assumption $\mathbf A$ and Cauchy-Schwarz's Inequality, 
\begin{align}
\esp\left | \int_0^{T\wedge \tau} \vartheta_{2}(s,Z_{s}){\rm d}X_{s} \right | & \le \esp\left | \int_0^{T\wedge\tau} \vartheta_{2}(s,Z_{s}){\rm d}\left(X_{s}-\Lambda_{Y,Z}(s)\right) \right | 
\nonumber\\
& \quad+ D^{2}\sqrt{T\mu([0,1]^{d})}\Vert \vartheta_{2}\Vert_{\mathbb L^{2}({\rm d}s\otimes\mu)}. 
\label{ref-th3-41}
\end{align}
Since $X-\Lambda_{Y,Z}$ is a martingale conditionally to $Y$ and $Z$, the Ito isometry \citep[see Theorem I.4.40 in][]{JaSh03} yields
\begin{align}
&\esp\left[\left( \int_0^{T\land\tau} \vartheta_{2}(s,Z_{s}){\rm d}\left(X_{s}-\Lambda_{Y,Z}(s)\right) \right)^{2} \,\Big |\, Y,Z\right]
\nonumber\\
&\quad=\esp\left[ \int_0^{T\land\tau} \vartheta^2_{2}(s,Z_{s}){\rm d}\langle X-\Lambda_{Y,Z}\rangle_s \,\Big |\, Y,Z\right],
\label{ref-th3-5}
\end{align}
where $\langle M\rangle$ stands for the predictable compensator of the martingale $M$. Observing that, conditionally on $Y,Z$, $X$ is a Poisson process with compensator $\Lambda_{Y,Z}$, we deduce that
$\langle X-\Lambda_{Y,Z}\rangle=\langle X\rangle=\Lambda_{Y,Z}$. As a result, 
\begin{align}
&\esp\left[ \int_0^{T\land\tau} \vartheta^2_{2}(s,Z_{s}){\rm d}\langle X-\Lambda_{Y,Z}\rangle_s \,\Big|\, Y,Z\right]
\nonumber\\
&\quad=\esp\left[ \int_0^{T\wedge\tau} \vartheta^2_{2}(s,Z_{s})\lambda_{Y}(s,Z_{s}){\rm d}s \,\Big |\, Y,Z\right].
\label{ref-th3-6}
\end{align}
Hence, 
\begin{equation}
\label{ref-th3-61}
\esp\left | \int_0^{T\wedge\tau} \vartheta_{2}(s,Z_{s}){\rm d}\left(X_{s}-\Lambda_{Y,Z}(s)\right) \right | \le D\Vert \vartheta_{2}\Vert_{\mathbb L^{2}({\rm d}s\otimes\mu)}.
\end{equation}
Combining \eqref{ref-th3-41} and \eqref{ref-th3-61} we deduce that 
\begin{equation}
\esp\left | \int_0^{T\wedge \tau} \vartheta_{2}(s,Z_{s}){\rm d}X_{s} \right |\leq D(1+D\sqrt{T\mu([0,1]^{d})})\| \vartheta_{2}\|_{\mathbb L^{2}({\rm d}s\otimes\mu)}.
\label{ref-th3-7}
\end{equation}
Putting together identities \eqref{ref-th3-0}-\eqref{ref-th3-3} and \eqref{ref-th3-7} yields
\begin{align*}
&\inf_{f\in\mathcal F_{B}}\left(A(f)-A(f^{\star})\right)\\
& \quad \leq 2D\sqrt{T\mu([0,1]^{d})}\| \vartheta_{1}\|_{\mathbb L^{2}({\rm d}s\otimes\mu)}+2D(1+D\sqrt{T\mu([0,1]^{d})})\| \vartheta_{2}\|_{\mathbb L^{2}({\rm d}s\otimes\mu)}\\
& \qquad +2\min_{\vert x\vert\leq B}\left\vert x-\ln\frac{p_{+}}{p_{-}}\right\vert,
\nonumber
\end{align*}
which concludes the proof by definition of $\vartheta_{1}$ and $\vartheta_{2}$. $\square$

\appendix
\section{Appendix}
\subsection{A general theorem for model selection}
The objective of this section is to recall a general model selection result due to \citet{BlBoMa03}.\\

Let $\mathcal X$ be a measurable space and let $\ell:\mathbb R\times\{-1,1\}\rightarrow\mathbb R$ be a loss function. Given a function $g:\mathcal X\rightarrow\mathbb R$, we let $\ell(g)$ be a shorthand notation for the function $(x,y)\in\mathbb R\times\{-1,1\}\mapsto \ell(g(x),y)$. Let $P$ be a probability distribution on $\mathcal X\times\{-1,1\}$ and let $\mathfrak G$ be a set of extended-real valued functions on $\mathcal X$ such that, for all $g\in\mathfrak G$, $\ell(g)\in\mathbb L^2(P)$. The target function $g^{\star}$ is defined as
$$g^{\star}\in\underset{g\in\mathfrak G}{\arg\min}\ P\ell(g).$$
Let $(\mathcal G_{k})_{k\ge1}$ be a countable family of models such that, for all $k\ge1$, $\mathcal G_{k}\subset \mathfrak G$. For each $k\ge 1$, we define the empirical risk minimizer $\hat g_{k}$ as
$$\hat g_{k}\in\underset{g\in\mathcal G_{k}}{\arg\min}\ P_{n}\ell(g).$$
If pen denotes a real-valued function on $\mathbb N^{\star}$, we let the penalized empirical risk minimizer $\hat g$ be defined by $\hat g_{\hat k}$, where
$$\hat k\in \underset{k\ge1}{\arg\min}\left[P_{n}\ell(\hat g_{k})+\mbox{pen$(k)$}\right].$$
Recall that a function $\mathbf d:\mathfrak G\times\mathfrak G\rightarrow\R_{+}$ is a pseudo-distance if $(i)$ $\mathbf d(g,g)=0$, $(ii)$ $\mathbf d(g,g')=\mathbf d(g',g)$, and $(iii)$ $\mathbf d(g,g')\le\mathbf d(g,g'')+\mathbf d(g'',g')$ for all $g$, $g'$, $g''$ in $\mathfrak G$. Also, a function $\psi:\mathbb R_{+}\rightarrow\mathbb R_{+}$ is said to be a sub-root function if $(i)$ it is nondecreasing and $(ii)$ the function $r\in\mathbb R_{+}\mapsto\psi(r)/\sqrt{r}$ is nonincreasing.
\begin{theo}[\citealp{BlBoMa03}]
\label{B}
Assume that there exist a pseu\-do-distance $\mathbf d$ on $\mathfrak G$, a sequence of sub-root functions $(\psi_{k})_{k\ge1}$, and two nondecreasing sequences $(b_{k})_{k\ge1}$ and $(C_{k})_{k\ge1}$ of real numbers such that
\begin{itemize}
\item[$(i)$] $\forall g,g'\in\mathfrak G:$ $P(\ell(g)-\ell(g'))^2\leq \mathbf d^2(g,g')$;
\item[$(ii)$] $\forall k\ge1$, $\forall g\in\mathcal G_{k}:$ $\mathbf d^2(g,g^{\star})\leq C_{k}P(\ell(g)-\ell(g^{\star}))$;
\item[$(iii)$] $\forall k\ge1$, $\forall g\in\mathcal G_{k}$, $\forall (x,y)\in\mathcal X\times\{-1,1\}:$ $ |\ell(g(x),y)|\leq b_{k}$;
\item[] and, if $r^{\star}_{k}$ denotes the solution of $\psi_{k}(r)=r/C_{k}$,
\item[$(iv)$] $\forall k\ge1$, $\forall g_{0}\in\mathcal G_{k}$, $\forall r\ge r^{\star}_{k}:$  
$$\esp\sup\left\{\left |(P_{n}-P)\left(\ell(g)-\ell(g_{0})\right)\right |: g\in\mathcal G_{k},\mathbf d^2(g,g_{0})\leq r\right\}\leq \psi_{k}(r).$$
\end{itemize}
Let $(x_{k})_{k\ge1}$ be a non\-increasing se\-quence such that $\sum_{k\ge 1}e^{-x_{k}}\leq 1$. Let $\delta>0$ and $K>1$ be two fixed real numbers. If \emph{pen$(k)$} denotes a penalty term satisfying
$$\forall k\ge1, \quad \emph{pen$(k)$}\ge 250K\frac{r^{\star}_{k}}{C_{k}}+\frac{\left(65KC_{k}+56b_{k}\right)\left(x_{k}+\delta+\ln 2\right)}{3n},$$
then, with probability at least $1-e^{-\delta}$, one has
$$P\left(\ell(\hat g)-\ell(g^{\star})\right)\le\frac{K+\frac{1}{5}}{K-1}\,\inf_{k\ge1}\,\left\{\inf_{g\in\mathcal G_{k}}P\left(\ell(g)-\ell(g^{\star})\right)+2\emph{pen$(k)$}\right\}.$$
\end{theo}

\subsection{Expected supremum of Rademacher processes}
Let $S$ be a measurable space and let $P$ be a probability measure on $S$. Let $\mathcal G$ be a class of functions $g:S\rightarrow \mathbb R$. The Rademacher process $(R_{n}(g))_{g\in\mathcal G}$ associated with $P$ and indexed by $\mathcal G$ is defined by
$$R_{n}(g)=\frac{1}{n}\sum_{i=1}^{n}\sigma_{i}g(Z_{i}),$$
where $\sigma_{1},\dots,\sigma_{n}$ are i.i.d. Rademacher ran\-dom vari\-ables, and $Z_{1},\dots,Z_{n}$ is a sequence of i.i.d. random variables, with distribution $P$ and independent of the $\sigma_{i}$'s. 
\medskip

We recall in this subsection a bound for the supremum of the Rademacher process defined by
$$\| R_{n}\|_{\mathcal G}=\sup_{g\in\mathcal G} \left | R_{n}(g)\right |,$$
which follows from the results of \citet{GiKo06}. Let $G$ be a measurable envelope for $\mathcal G$, i.e., a measurable function $G:S\rightarrow\mathbb R_{+}$ such that, for all $x\in S$, 
$$\sup_{g\in \mathcal G}\left | g(x)\right|\leq G(x).$$ 
Define $\| G\|=\sqrt{PG^2}$ and $\| G\|_{n}=\sqrt{P_{n}G^2}$, where $P_{n}=n^{-1}\sum_{i=1}^{n}\delta_{Z_{i}}$ stands for the empirical measure associated to $Z_{1},\dots,Z_{n}$. Finally, let $\sigma^2>0$ be a real number satisfying
$$\sup_{g\in\mathcal G}Pg^2\le\sigma^2\le\| G\|^2.$$

\begin{theo}[\citealp{GiKo06}]
\label{A}
Assume that the fun\-cti\-ons in $\mathcal G$ are uniformly bounded by a constant $U>0$. Assume, in addition, that there exist two constants $C$ and $V>0$ such that, for all $n\ge 1$ and all $0<\epsilon\le2$, 
$$N\left(\e\| G\|_{n},\mathcal G,\mathbb L^2\left(P_{n}\right)\right)\leq \left(\frac{C}{\e}\right)^V.$$ 
Then, for all $n\ge 1$, 
$$\esp\| R_{n}\|_{\mathcal G}\le\frac{c\sigma}{\sqrt n}\sqrt{V\ln\left(\frac{c'\| G\|}{\sigma}\right)}\lor \frac{8c^2UV}{n}\ln\left(\frac{c'\| G\|}{\sigma}\right)\lor \frac{cU}{9n}\sqrt{V\ln\left(\frac{c'\| G\|}{\sigma}\right)},$$
where $c=432$ and $c'=2e\lor C$. 
\end{theo}

\subsection{Some stochastic calculus results}

Up to the stopped part, the following result is a classical consequence of the Girsanov theorem. We give its proof for convenience of the reader. 

\begin{lem} 
\label{gigi} 
Let $\mu$ (resp., $\nu$) be the distribution of a Poisson process on $[0,T]$ with intensity $\lambda:[0,T]\to \mathbb R_+^\star$ (resp., with intensity 1) stopped after its $u$-th jump. Then, $\mu$ and $\nu$ are equivalent. Moreover, 
$$\nu({\rm d} x)=\exp \left(-\int_0^{T\wedge \tau(x)}\left(1-\lambda(s)\right){\rm d}s-\int_0^{T\wedge \tau(x)}\ln \lambda (s){\rm d}x_s\right)\mu({\rm d} x),$$
where, for all $x\in\mathcal X$, $\tau(x)=\inf\{t\in [0,T] : x_t=u\}$.
\end{lem}
\noindent {\sc Proof.} Consider the canonical  Poisson process $N=(N_t)_{t\in [0,T]}$ with intensity $\lambda$ on the filtered space $(\mathcal X,(\mathcal A_t)_{t\in [0,T]},\mathbb P)$, where $\mathcal A_{t}=\sigma(N_{s}: s\in[0,t])$, and let, for all $t\in [0,T]$,
$$\Lambda(t)=\int_0^t \lambda (s){\rm d}s \quad \mbox{and} \quad h(t)=\frac{1}{\lambda(t)}-1.$$
Recall that the process $M=(M_t)_{t\in [0,T]}$ defined by $M_t=N_t-\Lambda (t)$ is a martingale. The Dol\'eans-Dade exponential $\mathcal E=(\mathcal E_t)_{t\in [0,T]}$ of the martingale $h.M$ \citep[see, e.g., Theorem I.4.61 in][]{JaSh03} is defined for all $t\in [0,T]$ by 
\begin{align}
\mathcal E_t &=  e^{h.M_t}\prod_{s\leq t} (1+\Delta h.M_s)e^{-\Delta h.M_s}\nonumber \\
& =  \exp\left(-\int_0^t h(s)\lambda(s){\rm d}s+\int_0^t \ln\left(1+h(s)\right) {\rm d} N_s\right)\nonumber\\
& =  \exp\left(-\int_0^t \left(1-\lambda(s)\right){\rm d}s-\int_0^t \ln \lambda(s){\rm d}N_s\right),\label{doleans}
\end{align}
where $\Delta h.M_s=h.M_s-h.M_{s^-}=h.N_s-h.N_{s^-}$. Equivalently, $\mathcal E$ is the solution to the stochastic equation 
$$\mathcal E=1+{\mathcal E}^-.(h.M)=1+(\mathcal E^-h).M,$$
where $\mathcal E^-$ stands for the process defined by $\mathcal E^-_t=\mathcal E_{t^-}$. In particular, $\mathcal E$ is a martingale. Observe also, since $N$ is a counting process, that the quadratic covariation between $M$ and $\mathcal E$ is
$$[M,\mathcal E]=(\mathcal E^-h).[N,N]=(\mathcal E^-h).N.$$
Consequently,
$$[M,\mathcal E]-(\mathcal E^-h).\Lambda=(\mathcal E^-h).M$$
is a martingale. Since $(\mathcal E^-h).\Lambda$ is a continuous and adapted process, it is a predictable process and the predictable compensator of $[M,\mathcal E]$ takes the form $\langle M,\mathcal E\rangle = (\mathcal E^-h).\Lambda$. Now let $\mathbb Q$ be the measure defined by 
$${\rm d}\mathbb Q=\mathcal E_T {\rm d}\mathbb P.$$
Since the process $\mathcal E$ is a martingale, $\mathbb Q$ is a probability and in addition, for all $t\in [0,T]$, 
\begin{equation} \label{densityprocess} {\rm d}\mathbb Q_t=\mathcal E_t {\rm d}\mathbb P_t ,\end{equation}
where $\mathbb Q_t$ and $\mathbb P_t$ are the respective restrictions of $\mathbb Q$ and $\mathbb P$ to $\mathcal A_t$. It follows, by the Girsanov theorem \citep[see, e.g., Theorem III.3.11 in][]{JaSh03}, that the stochastic process $M-(\mathcal E^-)^{-1}.\langle M,\mathcal E\rangle$ is a $\mathbb Q$-martingale. But, for all $t\in [0,T]$,
\begin{align*}
M_t-(\mathcal E^-)^{-1}.\langle M,\mathcal E\rangle_t &= N_t-\Lambda(t)-h.\Lambda(t)\\
& =  N_t-(1+h).\Lambda(t)\\
& =  N_t-t.
\end{align*}
Therefore, the counting process $N$ is such that $(N_t-t)_{t\in [0,T]}$ is a $\mathbb Q$-martingale.  In consequence, by the Watanabe theorem \citep[see, e.g., Theorem IV.4.5 in][]{JaSh03}, this implies that  the distribution of $N$ under $\mathbb Q$ is that of a Poisson process with unit intensity. So, $\nu=\mathbb Q_{T\wedge \tau}$, where $\mathbb Q_{T\wedge \tau}$ is the restriction of $\mathbb Q$ to the stopped $\sigma$-field $\mathcal A_{T\wedge \tau}$. Moreover, by Theorem III.3.4 in \citet{JaSh03} and identity \eqref{densityprocess}, we have 
$${\rm d} \mathbb Q_{T\wedge \tau}=\mathcal E_{T\wedge \tau} {\rm d}\mathbb P_{T\wedge \tau},$$
where the definition of $\mathbb P_{T\wedge \tau}$ is clear. Since  $\mu=\mathbb P_{T\wedge \tau}$, the result is a consequence of identity \eqref{doleans}. $\Box$

\bibliography{biblio-bcpF}
\end{document}